\documentclass[12pt]{article}
\usepackage{mathrsfs}
       \usepackage{graphicx,amssymb,mathrsfs,amsmath,latexsym,amsfonts,titlesec,amscd,epsfig,cite,color,
                                                                                   }

\def\R{{\mathbb{R}}}

\def\BMO{{\mathrm{BMO}}}
\def\loc{{\mathrm{loc}}}

\DeclareMathOperator*{\essinf}{ess\, inf}
\DeclareMathOperator*{\esssup}{ess\, sup}

\newcommand{\upcite}[1]{\textsuperscript{\textsuperscript{\cite{#1}}}}
\numberwithin{equation}{section} 
\allowdisplaybreaks   
\usepackage{indentfirst} 
\usepackage[amsmath,thmmarks]{ntheorem} 

\theoremstyle{definition} 
\theoremheaderfont{\normalfont\rmfamily\bf}
\theorembodyfont{\normalfont\rm } \theoremindent0em
\theoremseparator{\hspace{-.5em}}
\newtheorem{theorem}{\indent
                  Theorem}[section]
    \newtheorem{lemma}{\indent  Lemma} [section]

\newtheorem{proposition}{\indent  Proposition}[section]

    \newtheorem{definition}{\indent  Definition} [section]
        \newtheorem{remark}{\indent  Remark}  

\newtheorem*{acknowledgments}{\indent Acknowledgments\quad }
\theoremstyle{nonumberplain} 
\theoremheaderfont{\bf\rmfamily} \theorembodyfont{\normalfont \rm }
\theoremsymbol{$\blacksquare$}
\newtheorem{proof}{\indent Proof}

\title{\bf\Large
Boundedness  for fractional Hardy-type operator  on  Herz-Morrey spaces with variable exponent }
\author{\textsc{Jianglong Wu } 
\\
\small {\it  Department of Mathematics, Mudanjiang Normal University, Mudanjiang, 157011, China}
\\
 \small { Bull. Korean Math. Soc. 51 (2014), No. 2, pp. 423-435} 
 }
\date{} 
\begin{document}

\maketitle

\footnote{ 
\textit{AMS} (2010) \textit{Mathematics Subject Classification}:
Primary 42B20; Secondary 47B38. }

\begin{abstract}
In this paper, the fractional Hardy-type operator of variable order $\beta(x)$ is shown to be bounded
from the  Herz-Morrey spaces $M\dot{K}_{p_{_{1}},q_{_{1}}(\cdot)}^{\alpha,\lambda}(\R^{n})$ with variable exponent $q_{1}(x)$ into the weighted space $M\dot{K}_{p_{_{2}},q_{_{2}}(\cdot)}^{\alpha,\lambda}(\R^{n},\omega)$, where $\omega=(1+|x|)^{-\gamma(x)}$ with some $\gamma(x)>0$ and $ 1/q_{_{1}}(x)-1/q_{_{2}}(x)=\beta(x)/n$ when $q_{_{1}}(x)$ is not necessarily constant at infinity. It is assumed that the exponent $q_{_{1}}(x)$ satisfies the logarithmic continuity condition both locally and at infinity that  $1<  q_{1}(\infty)\le q_{1}(x)\le( q_{1})_{+}<\infty~(x\in \R^{n})$.
\end{abstract}

\section{Introduction}

Let $f$ be a locally integrable function on $\R^{n}$. The $n$-dimensional Hardy operator is defined by
$$
\mathscr{H}(f)(x):= \frac{1}{|x|^{n}} \int_{|t|<|x|} f(t) \mathrm{d}t,\ \  \ x\in \R^{n}\setminus \{0\}.
$$

In 1995, Christ and Grafakos\upcite{CG} obtained the result for the
boundedness of $\mathscr{H}$ on $L^{p}(\R^{n})\ (1<p<\infty)$
spaces, and they also found the exact operator norms of $\mathscr{H}$ on
this space. In 2007, Fu et al\upcite{FLLW} gave the central $\BMO$ estimates for
commutators of  $n$-dimensional fractional and Hardy  operators. And recently, author\upcite{WW,W,WL,WZ,WZ1,ZW} also considers the boundedness for Hardy  operator and its commutator in (variable exponent) Herz-Morrey spaces.

The theory of variable exponent Lebesgue spaces is started by Orlicz (see \cite{O}, 1931) and Nakano (see \cite{N1,N2}, 1950 and 1951). In particular, the definition of Musielak-Orlicz spaces is clearly stated in \cite{N1}. However, the variable exponent function space, due to the failure of translation invariance and related properties, is very difficult to analyze.

Nowadays there is an evident increase of investigations related to both the theory of the spaces $L^{q(\cdot)}(\Omega)$ themselves and the operator theory in these spaces. This is caused by possible applications to models with non-standard local growth (in elasticity theory, fluid mechanics, differential equations, see for example \cite{R}, \cite{DR} and references therein) and is based on recent breakthrough result on boundedness of the Hardy-Littlewood maximal operator in these spaces.
 By virtue of the fine works\upcite{KR,CDF,CFMP,CFN,D1,D2,DHHMS,K,L,N,PR}, some
important conditions on variable exponent, for example, the $\log$-H\"{o}lder conditions and the
Muckenhoupt type condition, have been obtained.

Now, we define the n-dimensional  fractional Hardy-type operators of variable order $\beta(x)$  as follows.

\begin{definition}\ \   \label{def.1}
Let $f$ be a locally integrable function on $\R^{n},~ 0\le\beta(x)<n$. The n-dimensional  fractional Hardy-type operators of variable order $\beta(x)$ are defined by
 \begin{alignat}{2}
 \mathscr{H}_{\beta(\cdot)}(f)(x) &:= \frac{1}{|x|^{n-\beta(x)}} \int_{|t|<|x|} f(t)\mathrm{d}t,  \tag{\theequation a}  \label{equ.a}\\
\mathscr{H}^{\ast}_{\beta(\cdot)}(f)(x) &:=\int_{|t|\ge|x|} \frac{f(t)}{|t|^{n-\beta(x)}} \mathrm{d}t,  \tag{\theequation b}   \label{equ.b}
\end{alignat}
where $x\in \R^{n}\setminus \{0\}$.
\end{definition}

Obviously, when $\beta(x)=0$, $\mathscr{H}_{\beta(\cdot)}$ is just $\mathscr{H}$, and denote by $\mathscr{H}^{\ast}:=\mathscr{H}^{\ast}_{\beta(\cdot)}=\mathscr{H}^{\ast}_{0}$. And when $\beta(x)$ is constant, $\mathscr{H}_{\beta(\cdot)}$ and $\mathscr{H}^{\ast}_{\beta(\cdot)}$ will become $\mathscr{H}_{\beta}$ and $\mathscr{H}^{\ast}_{\beta}$\upcite{FLLW} respectively.

The Riesz-type potential operator of variable order $\beta(x)$ is defined by
\begin{equation} \label{equ.1}
I_{\beta(\cdot)}(f)(x)=\int_{\R^{n}} \frac{f(y)}{|x-y|^{n-\beta(x)}}  \mathrm{d}y,\  \ 0<\beta(x)<n.
\end{equation}
The boundedness of the operator $I_{\beta(\cdot)}$ from the space $L^{p(\cdot)}(\R^{n})$  with the variable exponent $p(x)$
into the space $L^{q(\cdot)}(\R^{n})$  with the limiting Sobolev exponent
 \begin{equation*}
\frac{1}{q(x)}=\frac{1}{p(x)}-\frac{\beta(x)}{n}
\end{equation*}
was an open problem for a long time. It was solved in the case of bounded
domains. First, in \cite{S}, in the case of bounded domains $\Omega$, there was proved a conditional result: the Sobolev theorem is valid for the potential operator  $I_{\beta(\cdot)}$ within the framework of the spaces $L^{p(\cdot)}(\Omega)$  with the variable exponent $p(x)$ satisfying the logarithmic Dini condition, if the maximal operator is bounded in the space $L^{p(\cdot)}(\Omega)$. In 2004,  Diening\upcite{D2} proved the boundedness of the maximal operator.   

In 2004, Diening\upcite{D} proved Sobolev's theorem for the potential  $I_{\beta}$ on
the whole space $\R^{n}$ assuming that $p(x)$ is constant at infinity ($p(x)$ is always constant outside some large ball) and satisfies the same logarithmic condition as in \cite{S}. Another progress for unbounded domains is the  result of Cruz-Uribe et al\upcite{CFN} on the boundedness of the maximal operator in unbounded domains for exponents $p(x)$ satisfying the logarithmic smoothness condition both locally and at infinity.

In \cite{KS}, Kokilashvili and Samko prove Sobolev-type theorem for the potential $I_{\beta(\cdot)}$ from the space $L^{p(\cdot)}(\R^{n})$  into the weighted space $L^{q(\cdot)}_{\omega}(\R^{n})$  with the power weight $\omega$ fixed to infinity, under the logarithmic condition for  $p(x)$ satisfied locally and at infinity, not supposing that  $p(x)$ is constant at infinity but assuming that $p(x)$ takes its minimal value at infinity.

 Motivated by the above results, we are to investigate mapping properties of the  fractional Hardy-type operators $\mathscr{H}_{\beta(\cdot)}$ and $\mathscr{H}^{\ast}_{\beta(\cdot)}$  within the framework of the Herz-Morrey spaces with variable exponent.

Throughout this paper, we will denote by $|S|$ the Lebesgue measure  and  by $\chi_{_{\scriptstyle S}}$ the characteristic function  for a measurable set  $S\subset\R^{n}$; $B(x,r)$ is the ball cenetered at $x$ and of radius $r$;~$B_{0}=B(0,1)$. $C$  denotes a constant that is independent of the main parameters involved but
whose value may differ from line to line. For any index $1< q(x)< \infty$, we denote by $q'(x)$ its conjugate index,
namely, $q'(x)=\frac{q(x)}{q(x)-1}$.  For $A\sim D$, we mean that there is a constant $C > 0$ such that$C^{-1}D\le A \le CD$.

\section{Preliminaries}

In this section, we give the definition of Lebesgue and Herz-Morrey spaces with variable exponent, and give basic properties
and useful lemmas.

\subsection{Function spaces with variable exponent}

Let $\Omega$ be a measurable set in $\R^{n}$ with $|\Omega|>0$. We first define Lebesgue spaces with variable exponent.

\begin{definition} \label{def.2}
\ \ Let ~$ q(\cdot): \Omega\to[1,\infty)$ be a measurable function.
 \renewcommand{\theenumi}{\Roman{enumi})}
\begin{enumerate}
 \item \ The Lebesgue spaces with variable exponent $L^{q(\cdot)}(\Omega)$ is defined by
  $$ L^{q(\cdot)}(\Omega)=\{f~ \mbox{is measurable function}:  F_{q}(f/\eta)<\infty ~\mbox{for some constant}~ \eta>0\}, $$
  where $F_{q}(f):=\int_{\Omega} |f(x)|^{q(x)} \mathrm{d}x$. The Lebesgue space $L^{q(\cdot)}(\Omega)$ is a Banach function space with respect to the norm
  \begin{equation*}
   \|f\|_{L^{q(\cdot)}(\Omega)}=\inf \Big\{ \eta>0:  F_{q}(f/\eta)=\int_{\Omega} \Big( \frac{|f(x)|}{\eta} \Big)^{q(x)} \mathrm{d}x \le 1 \Big\}.
\end{equation*}

  \item \ The space $L_{\loc}^{q(\cdot)}(\Omega)$ is defined by
  $$ L_{\loc}^{q(\cdot)}(\Omega)=\{f ~\mbox{is measurable}: f\in L^{q(\cdot)}(\Omega_{0}) ~\mbox{for all compact  subsets}~ \Omega_{0}\subset \Omega\}.  $$

 \item[III)]\  The weighted Lebesgue space $L_{\omega}^{q(\cdot)}(\Omega)$ is defined by as the set of all measurable functions for which $$\|f\|_{L^{q(\cdot)}_{\omega}(\Omega)}=\|\omega f\|_{L^{q(\cdot)}(\Omega)}<\infty.$$

    \end{enumerate} 
\end{definition}

Next we define some classes of variable  exponent functions. Given a function $f\in L_{\loc}^{1}(\R^{n})$, the Hardy-Littlewood maximal operator $M$ is defined by
$$Mf(x)=\sup_{r>0} r^{-n} \int_{B(x,r)} |f(y)| \mathrm{d}y,$$
where $B(x,r)=\{y\in \R^{n}: |x-y|<r\}$.

\begin{definition} \label{Def.3}\ \
 \ Given a measurable function $q(\cdot)$ defined on $\R^{n}$, we write
$$q_{-}:=\essinf_{x\in \R^{n}} q(x),\ \ q_{+}:= \esssup_{x\in \R^{n}} q(x).$$

\begin{list}{}{}
\item[(I)]\   $q'_{-}=\essinf\limits_{x\in \R^{n}} q'(x)=\frac{q_{+}}{q_{+}-1},\ \ q'_{+}= \esssup\limits_{x\in \R^{n}} q'(x)=\frac{q_{-}}{q_{-}-1}.$

\item[(II)]\ Denote by $\mathscr{P}(\R^{n})$ the set of all measurable functions $ q(\cdot): \R^{n}\to(1,\infty)$ such that
$$1< q_{-}\le q(x) \le q_{+}<\infty,\ \ x\in \R^{n}.$$

\item[(III)]\  The set $\mathscr{B}(\R^{n})$ consists of all  measurable functions  $q(\cdot)\in\mathscr{P}(\R^{n})$ satisfying that the Hardy-Littlewood maximal operator $M$ is bounded on $L^{q(\cdot)}(\R^{n})$.

\item[(IV)]\  The set $\mathscr{C}^{\log}_{0}(\R^{n})$ consists of all locally $\log$-H\"{o}lder continuous functions $ q(\cdot): \R^{n}\to(0,\infty)$ satisfies the condition
    \begin{equation}  \label{equ.4}
 |q(x)-q(y)| \le \frac{-C}{\ln(|x-y|)},  \ \ \  |x-y|\le 1/2,\ x,y \in \R^{n}.
\end{equation}

\item[(V)]\  The set $\mathscr{C}^{\log}_{\infty}(\R^{n})$ consists of all  $\log$-H\"{o}lder continuous at infinity functions $ q(\cdot): \R^{n}\to(0,\infty)$ satisfies the condition
\begin{equation}\label{equ.5}
 |q(x)-q(\infty)| \le \frac{C_{\infty}}{\ln(e+|x|)},  \ \ \ x \in \R^{n},
\end{equation}
where $q(\infty)=\lim_{|x|\to \infty}q(x)$.

\item[(VI)]\  Denote by  $\mathscr{C}^{\log}(\R^{n}):=\mathscr{C}^{\log}_{0}(\R^{n})\cap \mathscr{C}^{\log}_{\infty}(\R^{n})$ the set of all globally $\log$-H\"{o}lder continuous functions $ q(\cdot): \R^{n}\to(0,\infty)$.

\end{list}

\end{definition}

\begin{remark}\quad The $\mathscr{C}^{\log}_{\infty}(\R^{n})$ condition is equivalent to the uniform continuity condition
\begin{equation}\label{equ.6}
  |q(x)-q(y)| \le \frac{C}{\ln(e+|x|)},  \ \ \  |y|\ge|x|,\ x,y \in \R^{n}.
\end{equation}
The $\mathscr{C}^{\log}_{\infty}(\R^{n})$ condition was originally defined in this form in \cite{CFN}.
\end{remark}

  Next we define the Herz-Morrey spaces with  variable exponent. Let $B_{k}=B(0,2^{k})=\{x\in\R^{n}:|x|\leq 2^{k}\}, A_{k}=\ B_{k}\setminus  B_{k-1}$ and $\chi_{_{k}}=\chi_{_{A_{k}}}$ for $k\in \mathbb{Z}$.

\begin{definition} \label{def.4}\ \
Suppose that $\alpha\in \mathbb{R}, ~0\leq \lambda < \infty,~ 0<p< \infty$,
 $q(\cdot) \in \mathscr{P}(\mathbb{R}^{n})$. The Herz-Morrey space with  variable exponent $ M\dot{K}_{p, q(\cdot)}^{\alpha,
\lambda}(\mathbb{R}^{n} ) $ is definded by
  $$ M\dot{K}_{p, q(\cdot)}^{\alpha, \lambda}(\mathbb{R}^{n})=\Big\{f\in L_{\loc}^{q(\cdot)}(\mathbb{R}^{n}\backslash\{0\}):
 \|f\|_{M\dot{K}_{p, q(\cdot)}^{\alpha, \lambda}(\mathbb{R}^{n})}<\infty \Big\},  $$
where
  $$ \|f\|_{M\dot{K}_{p, q(\cdot)}^{\alpha, \lambda}(\mathbb{R}^{n})}=\sup_{k_{0}\in \mathbb{Z}}2^{-k_{0}\lambda}
 \Big(\sum_{k=-\infty}^{k_{0}}2^{k\alpha
 p}\|f\chi_{_{\scriptstyle k}}\|_{L^{^{q(\cdot)}}(\mathbb{R}^{n})}^{p} \Big)^{\frac{1}{p}}.  $$

\end{definition}

Compare the variable Herz-Morrey space $ M\dot{K}_{p, q(\cdot)}^{\alpha, \lambda}(\mathbb{R}^{n})$ with the variable Herz space  $\dot{K}_{q(\cdot)}^{\alpha,p}(\R^{n})$, where
$$\dot{K}_{q(\cdot)}^{\alpha,p}(\R^{n})= \Big\{f\in L_{\loc}^{q(\cdot)}(\mathbb{R}^{n}\backslash\{0\}):
\sum\limits_{k=-\infty}^{\infty}2^{k\alpha  p}\|f\chi_{_{k}}\|_{L^{q(\cdot)}(\R^{n})}^{p}<\infty \Big\},$$
Obviously, $M\dot{K}_{p,q(\cdot)}^{\alpha,0} (\R^{n})=\dot{K}_{q(\cdot)}^{\alpha,p}(\R^{n})$.

 In 2012,  Almeida and  Drihem\upcite{AD} discuss the boundedness of a wide class of sublinear operators on Herz spaces $K_{q(\cdot)}^{\alpha(\cdot),p}(\R^{n})$ and $\dot{K}_{q(\cdot)}^{\alpha(\cdot),p}(\R^{n})$ with variable exponent $\alpha(\cdot)$ and $q(\cdot)$.  In this paper, the author only considers Herz-Morrey space $ M\dot{K}_{p, q(\cdot)}^{\alpha(\cdot), \lambda}(\mathbb{R}^{n})$ with variable exponent $q(\cdot)$ but fixed $ \alpha\in \mathbb{R}$ and $p\in(0,\infty)$. However, for  the case of the exponent $\alpha(\cdot)$ is variable as well, which  can be found in the furthermore work for the author.

\subsection{Auxiliary propositions  and lemmas }

 In this part we   state some auxiliary propositions  and lemmas which will be needed for proving  our main theorems. And we only describe partial results  we need.

 \begin{proposition}\quad \label{Pro.1}
Let $ q(\cdot)\in \mathscr{P}(\R^{n})$.
 \begin{list}{}{}
 \item[(I)]\ If $ q(\cdot)\in \mathscr{C}^{\log}(\R^{n})$, 
 then we have $ q(\cdot)\in \mathscr{B}(\R^{n})$.
 \item[(II)]\  $q(\cdot)\in\mathscr{B}(\R^{n})$ if and only if $q'(\cdot)\in \mathscr{B}(\R^{n})$.
\end{list}
\end{proposition}

 The  first part in Proposition \ref{Pro.1} is independently due to Cruz-Uribe et al\upcite{CFN} and to Nekvinda\upcite{N}  respectively. The second of Proposition \ref{Pro.1}  belongs to  Diening\upcite{D1}~(see Theorem 8.1 or Theorem 1.2 in \cite{CFMP}).

\begin{remark}\quad \label{rem.2}
Since
$$|q'(x)-q'(y)|\le \frac{|q(x)-q(y)|}{(q_{-}-1)^{2}},$$
 it follows at once that if $q(\cdot)\in \mathscr{C}^{\log}(\R^{n})$, then so does $q'(\cdot)$---i.e., if the condition hold, then $M$ is bounded on
  $L^{q(\cdot)} (\R^{n})$ and $L^{q'(\cdot)} (\R^{n})$.  Furthermore, Diening has proved general results on Musielak-Orlicz spaces.

\end{remark}

The order $\beta(x)$ of the fractional Hardy-type operators in Definition \ref{def.1} is not assumed to be continuous. We assume that it is a measurable function on $\R^n$ satisfying the following assumptions
 \begin{equation} \label{equ.7}
\left.
 \begin{aligned}
         \beta_{0}:= \essinf_{x\in \R^{n}} \beta(x) & >0  \\
         \esssup_{x\in \R^{n}} p(x)\beta(x) &<n           \\
         \esssup_{x\in \R^{n}} p(\infty)\beta(x) &<n
                          \end{aligned} \right\}.
                          \end{equation}


In order to prove our main results, we need the Sobolev type theorem for the space $\R^{n}$ which was proved in ref. \cite{KS} for the
exponents $p(x)$ not necessarily constant in a neigbourhood of infinity, but with some ‘extra’
power weight fixed to infinity and under the assumption that $p(x)$ takes its minimal value at
infinity.
\begin{proposition}\ \ \label{pro.2}
Suppose that  $p(\cdot)\in \mathscr{C}^{\log}(\R^{n}) \cap\mathscr{P}(\R^{n})$. Let
\begin{equation}\label{equ.8}
1<  p(\infty)\le p(x)\le p_{+}<\infty,
\end{equation}
 and  $\beta(x)$ meet condition (\ref{equ.7}).
Then the following weighted Sobolev-type estimate is valid for the operator $I_{\beta(\cdot)}$:
 $$\Big\|(1+|x|)^{-\gamma(x)}I_{\beta(\cdot)}(f)\Big\|_{L^{q(\cdot)}(\R^{n})}\le C\|f\|_{L^{p(\cdot)}(\R^{n})},$$
where
\begin{equation*}
\frac{1}{q(x)}=\frac{1}{p(x)}-\frac{\beta(x)}{n}
\end{equation*}
is the Sobolev exponent and
\begin{equation}\label{equ.9}
\gamma(x)=C_{\infty} \beta(x) \Big(1-\frac{\beta(x)}{n}\Big) \le \frac{n}{4} C_{\infty},
\end{equation}
$C_{\infty}$ being the Dini-Lipschitz constant from (\ref{equ.5}) which $q(\cdot)$ is replaced by $p(\cdot)$.
\end{proposition}

\begin{remark}\quad \label{rem.3}
(i)\ \ If $\beta(x)$ satisfies the condition of type (\ref{equ.5}): $|\beta(x)-\beta(\infty)| \le \frac{C_{\infty}}{\ln(e+|x|)}$  $(x \in \R^{n})$,  then the weight   $(1+|x|)^{-\gamma(x)}$ is equivalent to the weight $(1+|x|)^{-\gamma(\infty)}$.

(ii)\ \ One can also treat operator (\ref{equ.1}) with $\beta(x)$ replaced by $\beta(y)$. In the
case of potentials over bounded domains $\Omega$ such potentials differ unessentially, if the function  $\beta(x)$ satisfies the smoothness logarithmic condition as (\ref{equ.4}), since
$$C_{1}|x-y|^{n-\beta(y)}\le |x-y|^{n-\beta(x)}\le C_{2}|x-y|^{n-\beta(y)}$$
in this case ( see \cite{S}, p. 277).

(iii)\ \ Under the assumptions of Proposition \ref{pro.2}, similar conclusion is also valid for the fractional maximal operator
$$M_{\beta(\cdot)}(f)(x)=\sup_{r>0}\frac{1}{|B(x,r)|^{n-\beta(x)}}\int_{B(x,r)} |f(y)| \mathrm{d}y.$$

(iv)\quad When $p(\cdot)\in \mathscr{P}(\R^{n})$,  the assumption that $p(\cdot)\in \mathscr{C}^{\log}(\R^{n})$ is equivalent to assuming  $1/p(\cdot)\in \mathscr{C}^{\log}(\R^{n})$, since
$$\Big|\frac{p(x)-p(y)}{(p_{+})^{2}}\Big| \le \Big|\frac{1}{p(x)}-\frac{1}{p(y)} \Big|=\Big|\frac{p(x)-p(y)}{p(x)p(y)}\Big|\le \Big|\frac{p(x)-p(y)}{(p_{-})^{2}}\Big|.$$
And further, $1/p(\cdot)\in \mathscr{C}^{\log}(\R^{n})$ implies that $1/q(\cdot)\in \mathscr{C}^{\log}(\R^{n})$ as well.
\end{remark}

The next lemma known as  the generalized H\"{o}lder's inequality on Lebesgue spaces with
variable exponent, and  the proof can be found in \cite{KR}.

\begin{lemma}\ (generalized H\"{o}lder's inequality) \label{Lem.1}
 \ \ Suppose that $ q(\cdot)\in \mathscr{P}(\R^{n})$, then for any $f\in L^{q(\cdot)}(\R^{n})$ and any $g\in L^{q'(\cdot)}(\R^{n})$, we have
 \begin{equation*}
 \int_{\R^{n}} |f(x)g(x)|\mathrm{d}x \le C_{q}\|f\|_{L^{q(\cdot)}(\R^{n})} \|g\|_{L^{q'(\cdot)}(\R^{n})},
\end{equation*}
where $C_{q}=1+1/q_{-}-1/q_{+}$.

\end{lemma}

The following  lemma can be found in \cite{I}. 

\begin{lemma}\ \ \label{Lem.2}
Let $ q(\cdot)\in \mathscr{B}(\R^{n})$.
\begin{list}{}{}
\item[(I)]\ \ Then there exist  positive constants $\delta\in (0,1)$ and $C>0$ such that
\begin{equation*}
 \frac{\|\chi_{S}\|_{L^{q(\cdot)}(\R^{n})}} {\|\chi_{B}\|_{L^{q(\cdot)}(\R^{n})}} \le C \bigg(\frac{|S|}{|B|} \bigg)^{\delta}
\end{equation*}
 for all balls $B$ in $\R^{n}$ and all measurable subsets $S\subset B$.

 \item[(II)]\ \ Then there exists a positive constant  $C>0$ such that
 \begin{equation*}  \label{equ.11}
 C^{-1}\le \frac{1}{|B|} \|\chi_{B}\|_{L^{q(\cdot)}(\R^{n})} \|\chi_{B}\|_{L^{q'(\cdot)}(\R^{n})} \le C
\end{equation*}
for all balls $B$ in $\R^{n}$.
 \end{list}
\end{lemma}

\begin{remark}\quad
(i)\ \  If $q_{_{1}}(\cdot),~q_{_{2}}(\cdot)\in \mathscr{C}^{\log}(\R^{n}) \cap \mathscr{P}(\R^{n})$, then we see that $ q'_{_{1}}(\cdot),~q_{_{2}}(\cdot) \in \mathscr{B}(\R^{n})$. Hence we can take positive constants $0<\delta_{1}<1/(q'_{_{1}})_{+},~0<\delta_{2}<1/(q_{_{2}})_{+}$  such that
\begin{equation}  \label{equ.12}
 \frac{\|\chi_{S}\|_{L^{q'_{1}(\cdot)}(\R^{n})}} {\|\chi_{B}\|_{L^{q'_{1}(\cdot)}(\R^{n})}} \le C \bigg(\frac{|S|}{|B|} \bigg)^{\delta_{1}},
 \ \ \  \frac{\|\chi_{S}\|_{L^{q_{2}(\cdot)}(\R^{n})}} {\|\chi_{B}\|_{L^{q_{2}(\cdot)}(\R^{n})}} \le C \bigg(\frac{|S|}{|B|} \bigg)^{\delta_{2}}
\end{equation}
 hold for all balls $B$ in $\R^{n}$ and all measurable subsets $S\subset B$ ( see \cite{WZ,I}).

(ii)\ \ On the other hand, Kopaliani\upcite{K} has proved the conclusion: If the exponent $q(\cdot)\in \mathscr{P}(\R^{n})$  equals to a constant outside some large ball, then $q(\cdot)\in \mathscr{B}(\R^{n})$ if and only if  $q(\cdot)$ satisfies the Muckenhoupt type condition
$$\sup_{Q:\hbox{cube}} \frac{1}{|Q|}\|\chi_{_{Q}}\|_{L^{q(\cdot)}(\R^{n})} \|\chi_{_{Q}}\|_{L^{q'(\cdot)}(\R^{n})} <\infty.$$
\end{remark}

\section{Main results and their proofs}

 Our main result can be stated as follows.

\begin{theorem}\quad\label{thm.1}
Suppose that $q_{_{1}}(\cdot)\in \mathscr{C}^{\log}(\R^{n}) \cap\mathscr{P}(\R^{n})$ satisfies condition (\ref{equ.8}), and $\beta(x)$ meet condition (\ref{equ.7}) which $p(\cdot)$ is replaced by $q_{_{1}}(\cdot)$. Define the variable exponent $q_{_{2}}(\cdot)$ by
$$\frac{1}{q_{_{2}}(x)}=\frac{1}{q_{_{1}}(x)}-\frac{\beta(x)}{n}.$$
Let $ 0<p_{_{1}}\le {p_{_{2}}} <\infty,~ \lambda\ge 0,~ \alpha<\lambda+ n\delta_{1}$,  where $\delta_{1}\in (0, 1/(q'_{1})_{+})$ is the constant appearing in (\ref{equ.12}).
Then  $$\Big\|(1+|x|)^{-\gamma(x)}\mathscr{H}_{\beta(\cdot)}(f) \Big\|_{M\dot{K}_{p_{_{2}},q_{_{2}}(\cdot)}^{\alpha,\lambda}(\R^{n})}\le C\|f\|_{M\dot{K}_{p_{_{1}},q_{_{1}}(\cdot)}^{\alpha,\lambda}(\R^{n})},$$
where $\gamma(x)$ is defined as in (\ref{equ.9}), and $C_{\infty}$ is the Dini-Lipschitz constant from (\ref{equ.5}) which $q_{_{1}}(\cdot)$ instead of $q(\cdot)$.

\end{theorem}

\begin{proof}\quad
For any $f\in {M\dot{K}_{p,  q(\cdot)}^{\alpha,\lambda}(\R^{n})}$, if we denote   $f_j:=f\cdot\chi_{j}=f\cdot\chi_{A_j}$ for each $j\in \mathbb{Z}$, then we can write
$$f(x)=\sum_{j=-\infty}^{\infty}f(x)\cdot\chi_{j}(x) =\sum_{j=-\infty}^{\infty}f_{j}(x).
$$

By (\ref{equ.a}) and Lemma \ref{Lem.1}, we have
\begin{eqnarray} \label{equ.3.1}
|\mathscr{H}_{\beta(\cdot)}(f)(x)\cdot\chi_{_{\scriptstyle k}}(x)| & \le& \frac{1}{|x|^{n-\beta(x)}} \int_{|t|<|x|} |f(t)|\mathrm{d}t \cdot\chi_{_{\scriptstyle k}}(x)    \le \frac{1}{|x|^{n-\beta(x)}} \int_{B_{k}} |f(t)|\mathrm{d}t \cdot\chi_{_{\scriptstyle k}}(x)   \nonumber\\
&\le& C 2^{^{-kn}}|x|^{\beta(x)} \Big(\sum_{j=-\infty}^{k}     \int_{A_{j}}|f(t)|\mathrm{d}t\Big)\cdot\chi_{_{\scriptstyle{k}}}(x)\\
& \le& C2^{^{-kn}} \sum_{j=-\infty}^{k} \|f_{j}\|_{L^{^{q_{_{1}}(\cdot)}}(\R^{n})} \|\chi_{_{\scriptstyle{j}}}\|_{L^{^{q'_{_{1}}(\cdot)}}(\R^{n})}  \cdot|x|^{\beta(x)}\chi_{_{\scriptstyle{k}}}(x). \nonumber
\end{eqnarray}

For Proposition \ref{pro.2}, we note that
\begin{equation} \label{equ.3.2}
\begin{split}
 I_{\beta(\cdot)}(\chi_{_{B_{k}}})(x)  &\ge I_{\beta(\cdot)}(\chi_{_{B_{k}}})(x)\cdot\chi_{_{B_{k}}}(x)    =   \int_{B_{k}}\frac{1}{|x-y|^{n-\beta(x)}} \mathrm{d}y  \cdot\chi_{_{B_{k}}}(x) \\
& \ge C|x|^{\beta(x)} \cdot\chi_{_{B_{k}}}(x)  \ge C|x|^{\beta(x)} \cdot\chi_{_{k}}(x).
\end{split}
\end{equation}

Using Proposition \ref{pro.2}, Lemma \ref{Lem.2},  (\ref{equ.12}),  (\ref{equ.3.1})  and (\ref{equ.3.2}), we have
\begin{equation} \label{equ.3.3}
\begin{split}
& \Big\|(1+|x|)^{-\gamma(x)}\mathscr{H}_{\beta(\cdot)}(f)\cdot\chi_{_{\scriptstyle k}}(\cdot)\Big\|_{L^{^{q_{_{2}}(\cdot)}}(\R^{n})} \\
  & \le C2^{^{-kn}} \sum_{j=-\infty}^{k} \|f_{j}\|_{L^{^{q_{_{1}}(\cdot)}}(\R^{n})} \|\chi_{_{\scriptstyle{j}}}\|_{L^{^{q'_{_{1}}(\cdot)}}(\R^{n})}\Big\|(1+|x|)^{-\gamma(x)}|\cdot|^{\beta(x)}\cdot\chi_{_{\scriptstyle k}}(\cdot)\Big\|_{L^{^{q_{_{2}}(\cdot)}}(\R^{n})}   \\
 & \le C2^{^{-kn}} \sum_{j=-\infty}^{k} \|f_{j}\|_{L^{^{q_{_{1}}(\cdot)}}(\R^{n})} \|\chi_{_{\scriptstyle{j}}}\|_{L^{^{q'_{_{1}}(\cdot)}}(\R^{n})}  \Big\|(1+|x|)^{-\gamma(x)}I_{\beta(\cdot)}(\chi_{_{B_{k}}})\Big\|_{L^{^{q_{_{2}}(\cdot)}}(\R^{n})}   \\
& \le C2^{^{-kn}} \sum_{j=-\infty}^{k}\|f_{j}\|_{L^{^{q_{_{1}}(\cdot)}}(\R^{n})} \|\chi_{_{\scriptstyle{j}}}\|_{L^{^{q'_{_{1}}(\cdot)}}(\R^{n})}  \|\chi_{_{B_{k}}}\|_{L^{^{q_{_{1}}(\cdot)}}(\R^{n})}   \\
& \le C2^{^{-kn}}\|\chi_{_{B_{k}}}\|_{L^{^{q_{_{1}}(\cdot)}}(\R^{n})} \sum_{j=-\infty}^{k}\|f_{j}\|_{L^{^{q_{_{1}}(\cdot)}}(\R^{n})} \|\chi_{_{B_{j}}}\|_{L^{^{q'_{_{1}}(\cdot)}}(\R^{n})}  \\
& \le C \sum_{j=-\infty}^{k}\|f_{j}\|_{L^{^{q_{_{1}}(\cdot)}}(\R^{n})} \frac{\|\chi_{_{B_{j}}}\|_{L^{^{q'_{_{1}}(\cdot)}}(\R^{n})}}{\|\chi_{_{B_{k}}}\|_{L^{^{q'_{_{1}}(\cdot)}}(\R^{n})}}    \le C\sum_{j=-\infty}^{k}2^{(j-k)n\delta_{1}} \|f_{j}\|_{L^{^{q_{_{1}}(\cdot)}}(\R^{n})}.
\end{split}
\end{equation}

Because of $0<p_{_{1}}/p_{_{2}}\le 1$, then we apply inequality
\begin{equation}     \label{equ.3.4}
\bigg(\sum_{i=-\infty}^{\infty}|a_{i}|\bigg)^{ p_{_{1}}/p_{_{2}}} \le \sum_{i=-\infty}^{\infty} |a_{i}|^{ p_{_{1}}/ p_{_{2}}},
\end{equation}
and obtain
\begin{eqnarray*}
&&\; \Big\|(1+|x|)^{-\gamma(x)}\mathscr{H}_{\beta(\cdot)}(f) \Big\|^{p_{_{1}}}_{M\dot{K}_{p_{_{2}},q_{_{2}}(\cdot)}^{\alpha,\lambda}(\R^{n})}  \\
&&\; = \sup_{k_{0}\in \mathbb{Z}}2^{-k_{0}\lambda {p_{_{1}}}}  \bigg(\sum_{k=-\infty}^{k_{0}}2^{k\alpha {p_{_{2}}}}  \Big\|(1+|x|)^{-\gamma(x)}\mathscr{H}_{\beta(\cdot)}(f)\cdot\chi_{_{\scriptstyle k}} (\cdot) \Big\|_{L^{q_{_{2}}(\cdot)}(\R^{n})}^{p_{_{2}}}\bigg)^{p_{_{1}}/p_{_{2}}}\\
&&\; \le  \sup_{k_{0}\in \mathbb{Z}}2^{-k_{0}\lambda p_{_{1}}}  \bigg(\sum_{k=-\infty}^{k_{0}}2^{k\alpha p_{_{1}}}
  \Big\|(1+|x|)^{-\gamma(x)}\mathscr{H}_{\beta(\cdot)}(f)\cdot\chi_{_{\scriptstyle k}}(\cdot) \Big\|_{L^{q_{_{2}}(\cdot)}(\R^{n})}^{p_{_{1}}}\bigg).
\end{eqnarray*}

On the other hand, note the following fact
\begin{eqnarray}
\begin{split} \label{equ.3.5}
\|f_{j}\|_{L^{^{q_{_{1}}(\cdot)}}(\R^{n})} &= 2^{-j\alpha}\Big(2^{j{\alpha}p_{_1}} \|f_{j}\|^{p_{_1}}_{L^{^{q_{_{1}}(\cdot)}}(\R^{n})}\Big)^{1/p_{_1}}\\
&\le  2^{-j\alpha}\bigg(\sum_{i=-\infty}^j 2^{i{\alpha}p_{_1}} \|f_{i}\|^{p_{_1}}_{L^{^{q_{_{1}}(\cdot)}}(\R^{n})}\bigg)^{1/p_{_1}}\\
&=  2^{j(\lambda-\alpha)}\bigg(2^{-j\lambda} \Big(\sum_{i=-\infty}^j 2^{i{\alpha}p_{_1}}
 \|f_{i}\|^{p_{_1}}_{L^{^{q_{_{1}}(\cdot)}} (\R^{n})}\Big)^{1/p_{_1}}\bigg)\\
&\le  C 2^{j(\lambda-\alpha)} \|f\|_{M\dot{K}_{p_{_{1}},q_{_{1}}(\cdot)}^{\alpha,\lambda}(\R^{n})}.
 \end{split}
\end{eqnarray}

Thus, combining (\ref{equ.3.3}) and (\ref{equ.3.5}), and using $\alpha< \lambda+n\delta_{1}$, it follows that
\begin{eqnarray*}
&&\; \Big\|(1+|x|)^{-\gamma(x)}\mathscr{H}_{\beta(\cdot)}(f) \Big\|^{p_{_{1}}}_{M\dot{K}_{p_{_{2}},q_{_{2}}(\cdot)}^{\alpha,\lambda}(\R^{n})}  \\
&&\; \le C\sup_{k_{0}\in \mathbb{Z}}2^{-k_{0}\lambda p_{_{1}}}  \bigg(\sum_{k=-\infty}^{k_{0}}2^{k\alpha p_{_{1}}}
  \Big(\sum_{j=-\infty}^{k} 2^{(j-k)n\delta_{1}} \|f_{j}\|_{L^{^{q_{_{1}}(\cdot)}}(\R^{n})} \Big)^{p_{_{1}}}\bigg) \\
&&\; \le C\|f\|^{p_{1}}_{M\dot{K}_{p_{_{1}},q_{_{1}}(\cdot)}^{\alpha,\lambda}(\R^{n})}\sup_{k_{0}\in \mathbb{Z}}2^{-k_{0}\lambda p_{_{1}}}  \bigg(\sum_{k=-\infty}^{k_{0}}2^{k\alpha p_{_{1}}}   \Big(\sum_{j=-\infty}^{k} 2^{(j-k)n\delta_{1}} 2^{j(\lambda-\alpha)} \Big)^{p_{_{1}}}\bigg) \\
&&\; \le C\|f\|^{p_{1}}_{M\dot{K}_{p_{_{1}},q_{_{1}}(\cdot)}^{\alpha,\lambda}(\R^{n})}\sup_{k_{0}\in \mathbb{Z}}2^{-k_{0}\lambda p_{_{1}}}  \bigg(\sum_{k=-\infty}^{k_{0}}2^{k\lambda p_{_{1}}}   \Big(\sum_{j=-\infty}^{k} 2^{(j-k)(n\delta_{1}+\lambda-\alpha)} \Big)^{p_{_{1}}}\bigg) \\
&&\; \le C\|f\|^{p_{1}}_{M\dot{K}_{p_{_{1}},q_{_{1}}(\cdot)}^{\alpha,\lambda}(\R^{n})} \sup_{k_{0}\in \mathbb{Z}}2^{-k_{0}\lambda p_{_{1}}}  \bigg(\sum_{k=-\infty}^{k_{0}}2^{k\lambda p_{_{1}}}\bigg)   \le C \|f\|^{p_{1}}_{M\dot{K}_{p_{_{1}},q_{_{1}}(\cdot)}^{\alpha,\lambda}(\R^{n})}.
\end{eqnarray*}

Consequently, the proof of Theorem \ref{thm.1} is completed.
\end{proof}


\begin{theorem}\quad\label{thm.2}
Suppose that $q_{_{1}}(\cdot)\in \mathscr{C}^{\log}(\R^{n}) \cap\mathscr{P}(\R^{n})$ satisfies condition (\ref{equ.8}), and $\beta(x)$ meet condition (\ref{equ.7}) which $q_{_{1}}(\cdot)$ instead of $p(\cdot)$. Define the variable exponent $q_{_{2}}(\cdot)$ by
$$\frac{1}{q_{_{2}}(x)}=\frac{1}{q_{_{1}}(x)}-\frac{\beta(x)}{n}.$$
Let $ 0<p_{_{1}}\le {p_{_{2}}} <\infty,~ \lambda\ge 0,~ \alpha>\lambda-n\delta_{2}$,  where $\delta_{2}\in (0, 1/(q_{2})_{+})$ is the constant appearing in (\ref{equ.12}).
Then  $$\Big\|(1+|x|)^{-\gamma(x)}\mathscr{H}^{\ast}_{\beta(\cdot)}(f) \Big\|_{M\dot{K}_{p_{_{2}},q_{_{2}}(\cdot)}^{\alpha,\lambda}(\R^{n})}\le C\|f\|_{M\dot{K}_{p_{_{1}},q_{_{1}}(\cdot)}^{\alpha,\lambda}(\R^{n})},$$
where $\gamma(x)$ is defined as in (\ref{equ.9}), and the Dini-Lipschitz constant from condition (\ref{equ.5}) which  $q(\cdot)$ is replaced by  $q_{_{1}}(\cdot)$.
\end{theorem}

\begin{proof}\quad
For simplicity, for any $f\in {M\dot{K}_{p,  q(\cdot)}^{\alpha,\lambda}(\R^{n})}$, we  write
$$f(x)=\sum_{j=-\infty}^{\infty}f(x)\cdot\chi_{j}(x) =\sum_{j=-\infty}^{\infty}f_{j}(x).
$$

By (\ref{equ.b}) and Lemma \ref{Lem.1}, we have
\begin{equation} \label{equ.3.6}
\begin{split}
&\Big|(1+|x|)^{-\gamma(x)}\mathscr{H}^{\ast}_{\beta(\cdot)}(f)(x)\cdot\chi_{_{\scriptstyle k}}(x)\Big| \le \int_{|t|\ge|x|} \frac{|f(t)|}{|t|^{n-\beta(x)}} \mathrm{d}t \cdot (1+|x|)^{-\gamma(x)}\chi_{_{\scriptstyle k}}(x)   \\
&\le C\int_{\R^{n}\setminus B_{k}} |f(t)| |x|^{\beta(x)-n} \mathrm{d}t \cdot (1+|x|)^{-\gamma(x)}\chi_{_{\scriptstyle k}}(x)   \\
& \le C\sum_{j=k+1}^{\infty} \int_{A_{j}} |f(t)| |x|^{\beta(x)-n} (1+|x|)^{-\gamma(x)}\mathrm{d}t \cdot \chi_{_{\scriptstyle k}}(x)   \\
& \le C\sum_{j=k+1}^{\infty}  \|f_{j}\|_{L^{^{q_{_{1}}(\cdot)}}(\R^{n})} \Big\|(1+|x|)^{-\gamma(x)} |\cdot|^{\beta(x)-n}\chi_{_{\scriptstyle{j}}} (\cdot) \Big\|_{L^{^{q'_{_{1}}(\cdot)}}(\R^{n})}  \cdot\chi_{_{\scriptstyle{k}}}(x).
\end{split}
\end{equation}
Similar to (\ref{equ.3.2}), we give
\begin{equation} \label{equ.3.7}
\begin{split}
 I_{\beta(\cdot)}(\chi_{_{B_{j}}})(x)  &\ge I_{\beta(\cdot)}(\chi_{_{B_{j}}})(x)\cdot\chi_{_{B_{j}}}(x)     =   \int_{B_{j}}\frac{1}{|x-y|^{n-\beta(x)}} \mathrm{d}y  \cdot\chi_{_{B_{j}}}(x)\\
& \ge C|x|^{\beta(x)} \cdot\chi_{_{B_{j}}}(x)  \ge C|x|^{\beta(x)} \cdot\chi_{_{j}}(x).
\end{split}
\end{equation}

 Since $q_{_{1}}(\cdot)\in \mathscr{C}^{\log}(\R^{n}) \cap\mathscr{P}(\R^{n})$  and $\beta(x)$ satisfy condition (\ref{equ.7}) and (\ref{equ.8})  which $q_{_{1}}(\cdot)$ instead of $p(\cdot)$. So,
applying Proposition \ref{pro.2},  Lemma \ref{Lem.2}, (\ref{equ.12}),  (\ref{equ.3.6})  and (\ref{equ.3.7}), we obtain
\begin{equation} \label{equ.3.8}
\begin{split}
& \Big\|(1+|x|)^{-\gamma(x)}\mathscr{H}^{\ast}_{\beta(\cdot)}(f)\cdot\chi_{_{\scriptstyle k}} (\cdot) \Big\|_{L^{^{q_{_{2}}(\cdot)}}(\R^{n})} \\
& \le C\sum_{j=k+1}^{\infty}  \|f_{j}\|_{L^{^{q_{_{1}}(\cdot)}}(\R^{n})} \|\chi_{_{\scriptstyle{k}}} \|_{L^{^{q_{_{2}}(\cdot)}}(\R^{n})}  \Big\|(1+|x|)^{-\gamma(x)} |\cdot|^{\beta(x)-n}\chi_{_{\scriptstyle{j}}} (\cdot)\Big\|_{L^{^{q'_{_{1}}(\cdot)}}(\R^{n})}    \\
  & \le C\sum_{j=k+1}^{\infty}   \|f_{j}\|_{L^{^{q_{_{1}}(\cdot)}}(\R^{n})} \|\chi_{_{\scriptstyle{k}}}\|_{L^{^{q_{_{2}}(\cdot)}}(\R^{n})} \cdot 2^{^{-jn}} \Big\|(1+|x|)^{-\gamma(x)} I_{\beta(\cdot)}(\chi_{_{B_{j}}})\Big\|_{L^{^{q'_{_{1}}(\cdot)}}(\R^{n})}   \\
& \le C\sum_{j=k+1}^{\infty} \|f_{j}\|_{L^{^{q_{_{1}}(\cdot)}}(\R^{n})} \|\chi_{_{B_{k}}}\|_{L^{^{q_{_{2}}(\cdot)}}(\R^{n})} \cdot 2^{^{-jn}} \|\chi_{_{B_{j}}}\|_{L^{^{q'_{_{2}}(\cdot)}}(\R^{n})}   \\
& \le C \sum_{j=k+1}^{\infty}\|f_{j}\|_{L^{^{q_{_{1}}(\cdot)}}(\R^{n})} \frac{\|\chi_{_{B_{k}}}\|_{L^{^{q_{_{2}}(\cdot)}}(\R^{n})}}{\|\chi_{_{B_{j}}}\|_{L^{^{q_{_{2}}(\cdot)}}(\R^{n})}}    \le C\sum_{j=k+1}^{\infty} 2^{(k-j)n\delta_{2}} \|f_{j}\|_{L^{^{q_{_{1}}(\cdot)}}(\R^{n})}.
\end{split}
\end{equation}

Because of $0<p_{_{1}}/p_{_{2}}\le 1$, therefore, applying (\ref{equ.3.4}), and combining  (\ref{equ.3.5}) and (\ref{equ.3.8}), and using $\alpha>\lambda-n\delta_{2}$, it follows that
\begin{eqnarray*}
&&\; \Big\|(1+|x|)^{-\gamma(x)}\mathscr{H}^{\ast}_{\beta(\cdot)}(f) \Big\|^{p_{_{1}}}_{M\dot{K}_{p_{_{2}},q_{_{2}}(\cdot)}^{\alpha,\lambda}(\R^{n})}  \\
&&\; \le C\sup_{k_{0}\in \mathbb{Z}}2^{-k_{0}\lambda p_{_{1}}}  \bigg(\sum_{k=-\infty}^{k_{0}}2^{k\alpha p_{_{1}}}
  \Big(\sum_{j=k+1}^{\infty} 2^{(k-j)n\delta_{2}} \|f_{j}\|_{L^{^{q_{_{1}}(\cdot)}}(\R^{n})} \Big)^{p_{_{1}}}\bigg) \\
&&\; \le C\|f\|^{p_{1}}_{M\dot{K}_{p_{_{1}},q_{_{1}}(\cdot)}^{\alpha,\lambda}(\R^{n})}\sup_{k_{0}\in \mathbb{Z}}2^{-k_{0}\lambda p_{_{1}}}  \bigg(\sum_{k=-\infty}^{k_{0}}2^{k\alpha p_{_{1}}}   \Big(\sum_{j=k+1}^{\infty} 2^{(k-j)n\delta_{2}}  2^{j(\lambda-\alpha)} \Big)^{p_{_{1}}}\bigg) \\
&&\; \le C\|f\|^{p_{1}}_{M\dot{K}_{p_{_{1}},q_{_{1}}(\cdot)}^{\alpha,\lambda}(\R^{n})}\sup_{k_{0}\in \mathbb{Z}}2^{-k_{0}\lambda p_{_{1}}}  \bigg(\sum_{k=-\infty}^{k_{0}}2^{k\lambda p_{_{1}}}   \Big(\sum_{j=k+1}^{\infty} 2^{(k-j)(n\delta_{2}+\alpha-\lambda)} \Big)^{p_{_{1}}}\bigg) \\
&&\; \le C\|f\|^{p_{1}}_{M\dot{K}_{p_{_{1}},q_{_{1}}(\cdot)}^{\alpha,\lambda}(\R^{n})} \sup_{k_{0}\in \mathbb{Z}}2^{-k_{0}\lambda p_{_{1}}}  \bigg(\sum_{k=-\infty}^{k_{0}}2^{k\lambda p_{_{1}}}\bigg)   \le C \|f\|^{p_{1}}_{M\dot{K}_{p_{_{1}},q_{_{1}}(\cdot)}^{\alpha,\lambda}(\R^{n})}.
\end{eqnarray*}

Consequently, the proof of Theorem \ref{thm.2} is completed.
\end{proof}


In particular, when $\gamma(x)=0$ and $\beta(\cdot)$ is constant exponent, the main results above are proved by Zhang and Wu in \cite{ZW}.

\begin{acknowledgments}
 The author cordially  thank the referees for their valuable suggestions and useful comments which have lead to the improvement of this paper. This work was partially supported   by the Fund  (No.12531720) of Heilongjiang Provincial Education Department, the Project (No.SY201313) of Mudanjiang Normal University   and NSF (No.11161042)  of China.
\end{acknowledgments}


\bigskip

\end{document}